\documentclass[a4j,12pt,reqno]{amsart}
\usepackage{amscd,epsfig,amsthm,graphicx}
\usepackage{amssymb}
\usepackage{amsmath}
\usepackage{amsthm}
\usepackage{amsfonts}
\usepackage[english]{babel}
\usepackage{epic,eepic}
\usepackage[flushleft]{paralist} 
\usepackage[utf8]{inputenc}
\usepackage{enumerate}
\usepackage{graphicx}
\usepackage[all]{xy} %xy-pic の使用

\makeindex			 %索引の作成
\newtheorem{theorem}{Theorem}
\newtheorem{lemma}[theorem]{Lemma}

\theoremstyle{definition}
			{Definition}

\newtheorem{remark}{Remark}

\makeatletter
    
    \@addtoreset{equation}{section}
\makeatother

\def\N{{\mathbb N}}
\def\Z{{\mathbb Z}} \def\Q{{\mathbb Q}}
\def\C{{\mathbb C}} 

\makeatletter
\def\@setcopyright{}
\def\serieslogo@{}
\makeatother

\begin{document}

% First we specify the top matter (author, title, etc).
%
% Note: All of the top matter items are optional and can be omitted.
% But you will probably want to specify at least the author and title
% and perhaps an abstract.

% author information

% first author 
\author{Takashi MIYAGAWA}      % 
%\address{Graduate School of Mathematics, Nagoya University, 
%		Chikusa-ku, Nagoya 464-8602, Japan}
%\email{d15001n@math.nagoya-u.ac.jp}

       % second author
       % \author{}
       % the address where the research was carried out
       % \address{}
       % current address, usually not needed because it is the same as the
       % regular address
       % \email{}
       % title

\title{Mean values of the Barnes double zeta-function}
%\renewcommand{\thefootnote}{\fnsymbol{footnote}}
%\begin{titlepage}
%\date{}
%\end{}

\begin{abstract}
In the study of order estimation of the Riemann zeta-function $ \zeta(s) = \sum_{n=1}^\infty n^{-s} $,
solving Lindel\"{o}f hypothesis is an important theme. 
As one of the relationships, asymptotic behavior of mean values has been studied.
Furthermore, the theory of the mean values is also noted in the double zeta-functions,
and the mean values of the Euler-Zagier type of double zeta-function and Mordell-Tornheim 
type of double zeta-function were studied. 
In this paper, we prove asymptotic formulas for mean square values of the Barnes double zeta-function 
 $ \zeta_2 (s, \alpha ; v, w ) = \sum_{m=0}^\infty \sum_{n=0}^\infty (\alpha+vm+wn)^{-s} $
 with respect to $ \mathrm{Im}(s) $ as $ \mathrm{Im}(s) \rightarrow + \infty $. 
\end{abstract}

\subjclass[2010]{Primary 11M32; Secondary 11B06}
\keywords{Barnes double zeta-function, Mean square values}
%\protect\footnotetext{2010{\it Mathematics Subject Classification.}\ \ Primary 11M32, Secondary 11M06. \\
%{\it Key words and phrases.}\ \ Barnes double zeta-function, Mean square values. \\
%{\it Address.}\ \ Graduate School of Mathematics, Nagoya University, Chikusa-ku, Nagoya 464-8602, Japan \\
%{\it E-mail.}\ \ d15001n@math.nagoya-u.ac.jp.}

\maketitle
% acknowledge support, etc
% \thanks{This research was partially supported by NSF grant DOA-123456789.}
% \thanks{We would like to thank our colleagues for their helpfulcriticism.}

% dedication
\dedicatory{}

% today's date, or fill in whatever date you prefer
\date{}

% This ends the top matter information.
% We can now tell LaTeX to display the top matter.

\maketitle

%\pagenumbering{roman}
%\tableofcontents
%\pagenumbering{arabic}

%%%%%%%% §1.Introduction  %%%%%%%%%%%%%%%%%%%%%%%%%%%%%%%%%%%%%%%%%%%%%%%%%
%%%%%%%%%%%%%%%%%%%%%%%%%%%%%%%%%%%%%%%%%%%%%%%%%%%%%%%%%%%%%%%%%%%%%%%%%%%%
%%%%%%%%%%%%%%%%%%%%%%%%%%%%%%%%%%%%%%%%%%%%%%%%%%%%%%%%%%%%%%%%%%%%%%%%%%%%
%%%%%%%%%%%%%%%%%%%%%%%%%%%%%%%%%%%%%%%%%%%%%%%%%%%%%%%%%%%%%%%%%%%%%%%%%%%%
%%%%%%%%%%%%%%%%%%%%%%%%%%%%%%%%%%%%%%%%%%%%%%%%%%%%%%%%%%%%%%%%%%%%%%%%%%%%
\section{Introduction and the statement of results}

The Barnes double zeta-function was first introduced by Barnes \cite{B1}
in the course of developing his theory of double gamma functions, and the double series of the 
form as
\begin{equation}\label{B_2-zeta}
	\zeta_2 (s,\alpha; v,w) = \sum_{m=0}^\infty \sum_{n=0}^\infty \frac{1}{(\alpha + vm + wn)^s}
\end{equation}
was introduced and studied in \cite{B2}.
As a subsequent research, multiple series of similar form as (\ref{B_2-zeta}) 
was introduced in connection with the theory of multiple gamma functions by Barnes \cite{B3}.

Let $ r $ be a positive integer, $ s = \sigma + it $ a complex variable, $ \alpha $
a real parameter, and $ w_j \ (j = 1, \ldots, r) $ complex parameters which are located
on one of the complex half-plane  divided by a straight line through the origin.
The Barnes multiple zeta-function $ \zeta_r(s, \alpha ; w_1, \ldots , w_r) $ is defined by
\begin{equation}\label{B-zeta}
       	\zeta_r (s, \alpha ; w_1, \ldots , w_r)
               = \sum_{m_1 = 0}^\infty \cdots \sum_{m_r = 0}^\infty
                 \frac{1}{(\alpha + w_1 m_1 + \cdots + w_r m_r)^s}
\end{equation}
where the series on the right-hand side is absolutely convergent for $ \mathrm{Re}(s) > r $,
and is continued meromorphically to the complex $ s $-plane, and its only singularities are
the simple poles located at $ s=j\ (j=1,2,\ldots,r) $.

In this paper, we focus on the case $ r=2 $ and $ (w_1,w_2) = (v,w) $ for any
$ v,w >0 $ of (\ref{B-zeta}), which is the Barnes double zeta-function (\ref{B_2-zeta})
and study the asymptotic behavior of
\[
	\int_1^T |\zeta_2(\sigma + it, \alpha;v,w)|^2 dt 
\]
as $ T \rightarrow +\infty $.

Let
\[
	\zeta_2^{[2]}(s_1,s_2, \alpha ; v,w) 
		= \mathop{\sum_{m_1, n_1, m_2, n_2 \geq 0}}
			\limits_{vm_1+wn_1 = vm_2+wn_2}
			\frac{1}{(\alpha + v m_1 + w n_1)^{s_1} 
							(\alpha + v m_2 + w n_2)^{s_2}},
\]
which is absolutely convergent for $ \mathrm{Re}(s_1 + s_2) > 2 $.
If $ v, w $ are linearly independent over $ \Q $, then $ vm_1 + wn_1 = vm_2 + wn_2 $
is equivalent to $ (m_1, n_1) = (m_2, n_2) $, and hence we have
\[
	\zeta_2^{[2]}(s_1,s_2, \alpha ; v,w) = \zeta_2(s_1+s_2, \alpha; v, w).
\]

\medskip

%================= Theorem1 ====================================================
\begin{theorem}\label{th:Main_Theorem1}
For $ s = \sigma + it \in \C $ with $ \sigma > 2 $, we have
\begin{equation*}
	\int_1^T | \zeta_2(s, \alpha ; v,w)|^2 dt 
	= \zeta_2^{[2]}(\sigma, \sigma, \alpha ; v,w)T + O(1)
\end{equation*}
as $ T \rightarrow +\infty $ .
\end{theorem}
%================= Theorem1 ====================================================

\medskip

%================= Theorem2 ====================================================
\begin{theorem}\label{th:Main_Theorem2}
For $ s = \sigma + it \in \C $ with $ 3/2 < \sigma \leq 2 $, we have
\begin{eqnarray*}
	&& \int_1^T | \zeta_2(s, \alpha ; v,w)|^2 dt 		\\
	&& \qquad \qquad \quad 
	= \zeta_2^{[2]}(\sigma, \sigma, \alpha ; v,w)T + 
	\begin{cases}
			O(T^{4-2\sigma} \log{T}) & 
				\left( 3/2 < \sigma \leq 7/4 \right)	\\
			O({T}^{1/2})	& (7/4 < \sigma \leq 2)
	\end{cases}
\end{eqnarray*}
as $ T \rightarrow +\infty $.
\end{theorem}
%================= Theorem2 ====================================================

\medskip

\begin{remark}
We mention here some recent results on mean values of double zeta-functions.
Matsumoto-Tsumura {\cite{MT}} treated the Euler double
zeta-function
\[
	\zeta_2 (s_1, s_2) = \sum_{m=1}^\infty \sum_{n=1}^\infty 
						\frac{1}{m^{s_1}(m_1 + n_1)^{s_2}}
\]
and gave some formulas which imply
\begin{equation}\label{Euler-Zagier}
	\int_2^T |\zeta_2 (\sigma_1 + it_1 , \sigma_2 + it_2)|^2 dt_2
	\sim \zeta_2^{[2]} (\sigma_1+it_1, 2\sigma_2) T	\qquad (T \rightarrow \infty)
\end{equation}
in some subsets in a region for $ \sigma_1 + \sigma_2 > 3/2 $, see {\cite{MT}} for detail. Here,
$ \zeta_2^{[2]} (\sigma_1+it_1, 2\sigma_2) $ is defined by
\[
	\zeta_2^{[2]}(s_1, s_2) = \sum_{k=2}^\infty 
			\left| \sum_{m=1}^{k-1}\frac{1}{m^{s_1}} \right|^2 \frac{1}{k^{s_2}}
\]
which is absolutely convergent for $ \mathrm{Re}(s_2)>1/2 $ and $ \mathrm{Re}(s_1+s_2)>3/2 $.
Ikeda-Matsuoka-Nagata {\cite{IMN}} extended the region of results of Matsumoto-Tsumura {\cite{MT}}, 
and further they gave some asymptotic formulas which imply
\[
	\int_2^T |\zeta_2 (\sigma_1 + it_1 , \sigma_2 + it_2)|^2 dt_2
	\asymp T\log{T}	\qquad (T \rightarrow \infty)
\]
on polygonal line
 $ \{(\sigma_1, \sigma_2)\, |\, \sigma_1 + \sigma_2 = 3/2 \ \mathrm{and} \ \sigma_2 > 1/2 \} 
\cup \{(\sigma_1, \sigma_2)\, |\, \sigma_1 > 1 \ \mathrm{and} \ \sigma_2 > 1/2 \} $.
Also, they gave similar results on
\[
	\int_2^T |\zeta_2 (\sigma + it, s_2)|^2 dt,
	\quad
	\int_2^T |\zeta_2 (\sigma_1 + it, \sigma_2 + it)|^2 dt,
\]
see {\cite{IMN}} for detail.
On the other hand, for the Mordell-Tornheim double zeta-function
\[
	\zeta_{MT,2}(s_1,s_2;s_3)
	= \sum_{m=1}^\infty \sum_{n=1}^\infty
	\frac{1}{m^{s_1}n^{s_2}(m+n)^{s_3}},
\]
Okamoto-Onozuka {\cite{OO}} obtained some results on the mean square values 
which imply
\begin{equation}\label{Mordell-Tornheim}
	\int_2^T |\zeta_{MT,2}(s_1,s_2; \sigma + it)|^2 dt
	\sim \zeta_{MT,2}^{[2]}(s_1, s_2; 2 \sigma) T \qquad (T \rightarrow \infty)
\end{equation}
in some subset in the region for $ \sigma_1 + \sigma_2 + \sigma > 3/2 $, here
$ \zeta_{MT,2}^{[2]}(s_1, s_2; 2 \sigma) $ is defined by
\[
	\zeta_{MT,2}^{[2]}(s_1, s_2; s) = \sum_{k=2}^\infty 
		\left| \sum_{m=1}^{k-1}\frac{1}{m^{s_1}(k-m)^{s_2}} \right|^2 \frac{1}{k^s},
\]
which is absolutely convergent for $ 2 \mathrm{Re}(s_1) + \mathrm{Re}(s) > 1,\,
 2 \mathrm{Re}(s_2) + \mathrm{Re}(s) > 1 $ and 
$ 2 \mathrm{Re}(s_1) + 2\mathrm{Re}(s_2) + \mathrm{Re}(s) > 3 $.
Theorem \ref{th:Main_Theorem1} and Theorem \ref{th:Main_Theorem2} are the
results corresponding to (\ref{Euler-Zagier}), (\ref{Mordell-Tornheim}) of the
Barnes double zeta-function (\ref{B_2-zeta}) version.
\end{remark}

%%%%%%%% §2.Proof of Theorem 1 %%%%%%%%%%%%%%%%%%%%%%%%%%%%%%%%%%%%%%%%%%%%%%
%%%%%%%%%%%%%%%%%%%%%%%%%%%%%%%%%%%%%%%%%%%%%%%%%%%%%%%%%%%%%%%%%%%%%%%%%%%%%%
%%%%%%%%%%%%%%%%%%%%%%%%%%%%%%%%%%%%%%%%%%%%%%%%%%%%%%%%%%%%%%%%%%%%%%%%%%%%%%
%%%%%%%%%%%%%%%%%%%%%%%%%%%%%%%%%%%%%%%%%%%%%%%%%%%%%%%%%%%%%%%%%%%%%%%%%%%%%%
%%%%%%%%%%%%%%%%%%%%%%%%%%%%%%%%%%%%%%%%%%%%%%%%%%%%%%%%%%%%%%%%%%%%%%%%%%%%%%
\section{Proof of Theorem 1}

In this section, we give a proof of Theorem \ref{th:Main_Theorem1}.

\medskip
\textbf{Proof of Theorem \ref{th:Main_Theorem1}}. 
Let $ \sigma + it \in \C $ with $ \sigma > 2 $. We first calculate
$ |\zeta_2 (s, \alpha; v,w)|^2 $. We have
\begin{eqnarray*}
	&& | \zeta_2(s, \alpha ; v,w)|^2
	  = \zeta_2(s, \alpha ; v,w) \overline{\zeta_2(s, \alpha ; v,w)}	\\
	&&\ 
		= \sum_{m_1=0}^\infty\sum_{n_1=0}^\infty
		\frac{1}{(\alpha + v m_1 + w n_1)^{\sigma + it}}
		\sum_{m_2=0}^\infty\sum_{n_2=0}^\infty
		\frac{1}{(\alpha + v m_2 + w n_2)^{\sigma - it}}	\\
	&&\ 
		= \mathop{\sum_{m_1, n_1, m_2, n_2 \geq 0}}
			\limits_{vm_1+wn_1 = vm_2+wn_2}
			\frac{1}{(\alpha + v m_1 + w n_1)^\sigma 
							(\alpha + v m_2 + w n_2)^\sigma}			\\
	&&\quad
		+ \mathop{\sum_{m_1, n_1, m_2, n_2 \geq 0}}
			\limits_{vm_1+wn_1 \neq vm_2+wn_2}
			\frac{1}{(\alpha + v m_1 + w n_1)^\sigma 
							(\alpha + v m_2 + w n_2)^\sigma}
			\left(\frac{\alpha + v m_2 + w n_2}
				{\alpha + v m_1 + w n_1}\right)^{it}		\\
	&&\ 
		= \zeta_2^{[2]}(\sigma,\sigma, \alpha ; v,w)			\\
	&&\quad
		+ \mathop{\sum_{m_1, n_1, m_2, n_2 \geq 0}}
			\limits_{vm_1+wn_1 \neq vm_2+wn_2}
			\frac{1}{(\alpha + v m_1 + w n_1)^\sigma 
							(\alpha + v m_2 + w n_2)^\sigma}
			\left(\frac{\alpha + v m_2 + w n_2}
				{\alpha + v m_1 + w n_1}\right)^{it}.	
\end{eqnarray*}
Hence we have
\begin{eqnarray*}
	&& \int_1^T | \zeta_2(s, \alpha ; v,w)|^2 dt = \zeta_2^{[2]}(\sigma,\sigma, \alpha ; v,w)(T-1)	\\
	&& \qquad \quad
		+ \mathop{\sum_{m_1, n_1, m_2, n_2 \geq 0}}
			\limits_{vm_1+wn_1 \neq vm_2+wn_2}
			\frac{1}{(\alpha + v m_1 + w n_1)^\sigma(\alpha + v m_2 + w n_2)^\sigma}	\\
	&& \qquad \qquad
			\times
			\frac{e^{iT\log\{(\alpha + v m_2 + w n_2)/(\alpha + v m_1 + w n_1)\}}
				- e^{i\log\{(\alpha + v m_2 + w n_2)/(\alpha + v m_1 + w n_1)\}}}
			{i \log\{(\alpha + v m_2 + w n_2)/(\alpha + v m_1 + w n_1)\}}.
\end{eqnarray*}
The second term on the right-hand side is
\begin{eqnarray*}
	&\ll & \mathop{\sum_{m_1, n_1, m_2, n_2 \geq 0}}
			\limits_{vm_1+wn_1 < vm_2+wn_2}
			\frac{1}
				{(\alpha + v m_1 + w n_1)^\sigma(\alpha + v m_2 + w n_2)^\sigma}	\\
	&& \qquad \qquad \qquad
		\times \frac{1}{\log\{(\alpha + v m_2 + w n_2)/(\alpha + v m_1 + w n_1)\}} 	\\
	&= & \left( \mathop{\sum_{m_1, n_1, m_2, n_2 \geq 0}}
			\limits_{\alpha + vm_1 + wn_1 < \alpha + vm_2 + wn_2 < 2(\alpha + vm_1 + wn_1)}
		+ \mathop{\sum_{m_1, n_1, m_2, n_2 \geq 0}}
			\limits_{ \alpha + vm_2 + wn_2  \geq 2(\alpha + vm_1 + wn_1)} 	
		\right) \\
	&& \quad 
		\frac{1}{(\alpha + v m_1 + w n_1)^\sigma(\alpha + v m_2 + w n_2)^\sigma
				\log\{(\alpha + v m_2 + w n_2)/(\alpha + v m_1 + w n_1)\} }.
\end{eqnarray*}
We denote the right-hand side by $ V_1 + V_2 $. Then we have
\begin{eqnarray*}
	V_2 &\ll& \mathop{\sum_{m_1, n_1, m_2, n_2 \geq 0}}
			  \limits_{\alpha + vm_2 + wn_2  \geq 2(\alpha + vm_1 + wn_1) }
			  \frac{1}
				{(\alpha + vm_1 + wn_1)^\sigma (\alpha + vm_2 + wn_2)^\sigma}	\\
		& \leq & \sum_{m_1 \geq 0} \sum_{n_1 \geq 0} \frac{1}{(\alpha+vm_1+wn_1)^\sigma}
			 \sum_{m_2 \geq 0} \sum_{n_2 \geq 0} \frac{1}{(\alpha+vm_2+wn_2)^\sigma}
		= O(1).
\end{eqnarray*}
Next we consider the order of $ V_1 $. The range of $ n_2 $ satisfying the inequalities 
$ \alpha + vm_1 + wn_1 < \alpha + vm_2 + wn_2 < 2(\alpha + vm_1 + wn_1) $
of the condition on the sum $ V_1 $ is
\[
	\frac{v}{w}(m_1 - m_2) + n_1 < n_2 < \frac{\alpha}{w} + \frac{v}{w}(2m_1-m_2) + 2n_1.
\]
Let $ \varepsilon = \varepsilon(m_1,m_2,n_1),  \delta = \delta (m_1,m_2,n_1) $ be
the quantities satisfying $ 0 \leq \varepsilon, \delta < 1 $ and
$ (v/w)(m_1-m_2) + n_1 + \varepsilon \in \Z $ and 
$ \alpha/w + (v/w)(2m_1-m_2) + 2n_1 - \delta \in \Z $.	 
Then $ K = \alpha/w + (v/w)m_1 + n_1 - \varepsilon - \delta $ is an integer, 
and $ n_2 $ can be rewritten as
\[
	n_2 = \frac{v}{w}(m_1 - m_2) + n_1 + \varepsilon + k \quad 
	(\mathrm{for\ some\ } k=0,\ 1,\ 2,\ \ldots,\ K).
\] 
Since $ K \asymp \alpha + vm_1 + wn_1 \asymp 1+m_1 + n_1 $, we have
\begin{eqnarray*}
	\log{\frac{\alpha + v m_2 + w n_2}{\alpha + v m_1 + w n_1}}
	%&=& \log{\frac{\alpha + v m_1 + w n_1 + wk + w \varepsilon }{\alpha + v m_1 + w n_1}}
	= 	\log\left(1 + \frac{w \varepsilon + wk}{\alpha + v m_1 + w n_1}\right)
	\asymp  \frac{w \varepsilon + wk}{\alpha + v m_1 + w n_1},
\end{eqnarray*}
and hence
\begin{eqnarray*}
	V_1 &\ll & \sum_{m_1 \geq 0} \sum_{n_1 \geq 0} 
			\sum_{0 \leq m_2 \ll K} \sum_{0 \leq k \leq K}
			\frac{1}{(\alpha + vm_1 +wn_1)^\sigma}	\\
		&&	\qquad \qquad \times \frac{1}{(\alpha + vm_1 +wn_1 + wk + w \varepsilon )^\sigma}
						\times \frac{\alpha + v m_1 + w n_1}{wk + w \varepsilon }		\\
		&\ll& \sum_{m_1 \geq 0} \sum_{n_1 \geq 0} \sum_{0 \leq m_2 \ll K} 
			\frac{\log{K}}{(\alpha + vm_1 + wn_1)^{2\sigma -1}}							\\
		&\ll& \sum_{m_1 \geq 0} \sum_{n_1 \geq 0} 
			\frac{\log{(\alpha + vm_1 + wn_1)}}{(\alpha + vm_1 + wn_1)^{2\sigma -2}}	\\
		&\ll& \sum_{m=0}^\infty \sum_{n=0}^\infty \frac{\log{(2+m+n)}}{(1+m+n)^{2\sigma -2}}
		\ll 1,
\end{eqnarray*}
provided that $ \sigma > 2 $. Therefore the proof of Theorem \ref{th:Main_Theorem1} is complete.
\qed

\bigskip

%%%%%%%% §3. The approximation theorem %%%%%%%%%%%
%%%%%%%%%%%%%%%%%%%%%%%%%%%%%%%%%%%%%%%%%%%%%%%%%%%%%%%%%%%%%%%%%%%%%%%%%%%%
%%%%%%%%%%%%%%%%%%%%%%%%%%%%%%%%%%%%%%%%%%%%%%%%%%%%%%%%%%%%%%%%%%%%%%%%%%%%
%%%%%%%%%%%%%%%%%%%%%%%%%%%%%%%%%%%%%%%%%%%%%%%%%%%%%%%%%%%%%%%%%%%%%%%%%%%%
\section{The approximation theorem}

Let $ \sigma_1 > 0,\ x \geq 1 $ and $ C > 1 $. Suppose $ s = \sigma + it \in \C $
with $ \sigma \geq \sigma_1 $ and $ |t| \leq 2\pi x/C $. Then
\begin{equation}\label{order_of_zeta}
	\zeta(s) = \sum_{1 \leq n \leq x} \frac{1}{n^s} 
					- \frac{x^{1-s}}{1-s} + O(x^{-\sigma})
					\qquad (x \rightarrow \infty).
\end{equation}
This asymptotic formula has been proved by Hardy and Littlewood 
(see Theorem 4.11 in Titchmarsh \cite{Tit}).
Here we prove an analogue of (\ref{order_of_zeta}) for the case of the Barnes 
double zeta-functions as follows.

%================== Thm3 =========================================
\begin{theorem}\label{order_of_Barnes_zeta}
Let $ 1 < \sigma_1 < \sigma_2,\ x \geq 1 $ and $ C > 1 $.
Suppose $ s = \sigma + it \in \C $ with $ \sigma_1 < \sigma < \sigma_2 $
and $ |t| \leq 2\pi x/C $. Then
\begin{eqnarray}
	&& \zeta_2(s,\alpha;v,w) = \sum_{0 \leq m \leq x} \sum_{0 \leq n \leq x} 
				\frac{1}{(\alpha + vm + wn)^s}		\nonumber \\
	&& \qquad\qquad	+ \frac{(\alpha + vx)^{2-s} + (\alpha + wx)^{2-s} - (\alpha + vx + wx)^{2-s}}
				{vw(s-1)(s-2)}  + O(x^{1-\sigma} )	\label{app_of_Barnes_zeta}  
\end{eqnarray}
as $ x \rightarrow \infty $. 
\end{theorem}
%====================================================================

%================== Lemmma4 =========================================
\begin{lemma}[Lemma 4.10 in \cite{Tit}]\label{exp sum}
Let $ f(\xi) $ be a real function with a continuous and steadily decreasing derivative
$ f'(\xi) $ in $ (a, b) $, and let $ f'(b) = c, \ f'(a) = d $. 
Let $ g(\xi) $ be a real positive decreasing function with a continuous 
derivative $ g'(\xi) $, satisfying that $ | g'(\xi) | $  is steadily decreasing.
Then
\begin{eqnarray}\label{exponential-sum}
	\sum_{a < n \leq b} g(n) e^{2\pi i f(n)}
	&=& \mathop{\sum_{\nu \in \Z}}\limits_{c-\varepsilon < \nu < d + \varepsilon}
	  \int_a^b g(\xi) e^{2\pi i (f(\xi) - \nu \xi)} d\xi	\nonumber \\
	&& \quad + O(g(a) \log(d - c +2) ) + O(|g'(a)|)
\end{eqnarray}
for an arbitrary $ \varepsilon \in (0, 1) $.
\end{lemma}
%=======================================================================
\textbf{Proof of Theorem \ref{order_of_Barnes_zeta}}.
Let $ N \in \N $ be sufficiently large. Then we have
\begin{eqnarray*}
	&& \sum_{m=0}^\infty \sum_{n=0}^\infty 
			\frac{1}{(\alpha + v m + w n)^s}	\\			
	&& \qquad \quad = \left( \sum_{m=0}^N \sum_{n=0}^N 
				+ \sum_{m=0}^N \sum_{n=N+1}^\infty + \sum_{m=N+1}^\infty \sum_{n=0}^N 
				+ \sum_{m=N+1}^\infty \sum_{n=N+1}^\infty \right)
				\frac{1}{(\alpha + v m + w n)^s}.
\end{eqnarray*}
We denote the second, the third and the fourth term on the right-hand side by 
$ A_1, A_2 $ and $ A_3 $, respectively. By the Euler-Maclaurin summation formula
(see Equation (2.1.2) in \cite{Tit}), we have for any $ a, b \in \Z $ with
$ 0 < a < b $,
\begin{eqnarray*}
&& \sum_{m=a+1}^b \frac{1}{(\alpha + v m + w n)^s}
	= \frac{(\alpha + v b + w n)^{1-s}-(\alpha + v a + w n)^{1-s}}{v(1-s)}	\\
&& \qquad \qquad -vs \int_a^b \frac{x-[x]-1/2}{(\alpha + v x + w n)^{s+1}} dx		
	+\frac{1}{2} 	
		\left\{ (\alpha + v b + w n)^{-s}-(\alpha + v a + w n)^{-s} \right\},	\\
&& \sum_{n=a+1}^b \frac{1}{(\alpha + v m + w n)^s}
	= \frac{(\alpha + v m + w b)^{1-s}-(\alpha + v m + w a)^{1-s}}{w(1-s)}	\\
&& \qquad \qquad -ws \int_a^b \frac{x-[x]-1/2}{(\alpha + v m + w x)^{s+1}} dx		
	+\frac{1}{2} 
		\left\{ (\alpha + v m + w b)^{-s}-(\alpha + v n + w a)^{-s} \right\}.
\end{eqnarray*}
If we take $ a = N $ and let $ b \rightarrow \infty $, we have
\begin{eqnarray*}
	&&\sum_{n=N+1}^\infty \frac{1}{(\alpha + vm + wn)^s}	\\
	&& = \frac{1}{w(s-1)}\cdot \frac{1}{(\alpha + vm + wN)^{s-1}}	\\
	&& \qquad \qquad \qquad
		- ws \int_N^\infty \frac{x-[x]-1/2}{(\alpha + vm + wx)^{s+1}}dx	
			- \frac{1}{2} \frac{1}{(\alpha + vm +wN)^s}		\\
	&& = \frac{1}{w(s-1)}\cdot \frac{1}{(\alpha + vm + wN)^{s-1}} -
		\frac{1}{2} \frac{1}{(\alpha + vm +wN)^s}
			+ O \left(N^{-\sigma} \right),
\end{eqnarray*}
for $ \sigma >1 $, uniformly in $ m = 0,1, \ldots $. Therefore we have
\begin{eqnarray*}
	A_1 &=& \sum_{m=0}^N \sum_{n=N+1}^\infty \frac{1}{(\alpha + v m + w n)^s}	\\
		&=& \sum_{m=0}^N \left\{ \frac{1}{w(s-1)}\cdot \frac{1}{(\alpha + vm + wN)^{s-1}}
				- \frac{1}{2} \cdot \frac{1}{(\alpha + vm +wN)^s}
			+ O \left(N^{-\sigma} \right) \right\}	\\
		&=& \frac{1}{w(s-1)}\sum_{m=1}^N \frac{1}{(\alpha + vm + wN)^{s-1}} -
			\frac{1}{2} \sum_{m=1}^N \frac{1}{(\alpha + vm +wN)^s} \\
		& & \qquad \qquad + \frac{1}{w(s-1)} \cdot \frac{1}{(\alpha + wN)^{s-1}} 
			- \frac{1}{2} \cdot \frac{1}{(\alpha + wN)^s} + O(N^{1- \sigma}).
\end{eqnarray*}
Applying again the formula (2.1.2) in \cite{Tit} to the first term and the second term on the 
right-hand side of the above, we obtain
\begin{eqnarray*}
	A_1 &=& \frac{1}{vw(s-1)(s-2)} 
		\left\{ \frac{1}{(\alpha + wN)^{s-2}} - \frac{1}{(\alpha +vN+wN)^{s-2}} \right\}	\\
		&&  - \frac{v}{w} \int_0^N \frac{x-[x]-1/2}{(\alpha + vx + wN)^s}dx 	\\
		&&	- \frac{1}{2w(s-1)} 
			\left\{  \frac{1}{(\alpha + wN)^{s-1}} - \frac{1}{(\alpha +vN+wN)^{s-1}} \right\}	\\
		&& - \frac{1}{2v(s-1)} 
			\left\{\frac{1}{(\alpha + wN)^{s-1}} - \frac{1}{(\alpha +vN+wN)^{s-1}} \right\} \\
		&& + \frac{vs}{2} \int_0^N \frac{x-[x]-1/2}{(\alpha + vx + wN)^{s+1}}dx 
			+ \frac{1}{4} \left\{  \frac{1}{(\alpha + wN)^s} - \frac{1}{(\alpha +vN+wN)^s} \right\}	\\
		&& + \frac{1}{w(s-1)} \cdot \frac{1}{(\alpha + wN)^{s-1}} 
			- \frac{1}{2} \cdot \frac{1}{(\alpha + wN)^s} + O(N^{1-\sigma})	\\
		&=& \frac{(\alpha + wN)^{2-s} - (\alpha + vN + wN)^{2-s}}{vw(s-1)(s-2)} + O(N^{1-\sigma}).
\end{eqnarray*}
Applying the same method to $ A_2 $ and $ A_3 $, we obtain
\begin{eqnarray*}
	A_2 &=& \frac{(\alpha + vN)^{2-s} - (\alpha + vN + wN)^{2-s}}{vw(s-1)(s-2)} + O(N^{1-\sigma})
	\qquad (\sigma > 1)	\\
	A_3 &=& \frac{(\alpha + vN + wN)^{2-s}}{vw(s-1)(s-2)} + O(N^{1-\sigma}) \qquad (\sigma > 1).
\end{eqnarray*}
Therefore we have
\begin{eqnarray}
	&& \zeta_2(s,\alpha;v,w) =  \sum_{m=0}^N \sum_{n=0}^N \frac{1}{(\alpha + v m + w n)^s}	
		\nonumber \\
	&& \qquad \qquad
	+ \frac{(\alpha + vN)^{2-s} + (\alpha + wN)^{2-s} - (\alpha + vN + wN)^{2-s}}{vw(s-1)(s-2)} 
	+ O(N^{1-\sigma}) 	\label{eq:sum_{x<N}}
\end{eqnarray}
for $ \sigma > 1 $. Next we consider the double sum on the right-hand side of (\ref{eq:sum_{x<N}}).
First we divide the sum as follows:
\begin{eqnarray*}
	&&\sum_{m=0}^N \sum_{n=0}^N 
			\frac{1}{(\alpha + v m + w n)^s}	\\	
	&& = \left( \sum_{m \leq x} \sum_{n \leq x} 
				+ \sum_{m \leq x} \sum_{x < n \leq N} + \sum_{n \leq x} \sum_{x < m \leq N} 
				+ \sum_{x < m \leq N} \sum_{x < n \leq N} \right)
				\frac{1}{(\alpha + v m + w n)^s}.
\end{eqnarray*}				
We denote the second, the third and the fourth term on the right-hand side by $ B_1 $, $ B_2 $ 
and $ B_3 $, respectively. Fix $ m \in \N $, set
\[
	f(\xi) = \frac{t}{2\pi}\log(\alpha + vm + w\xi), \ g(\xi) = (\alpha + vm + w\xi)^{-\sigma}
\]
and take $ (a,b) = (x,N) $ in Lemma \ref{exp sum}. Then we have
\[
	(c, d) = \left( \frac{tw}{2\pi(\alpha + vm + wN)}, \frac{tw}{2\pi(\alpha + vm + wx)} \right).
\]
We see that
\[
	|f'(x)| = \left| \frac{tw}{2\pi(\alpha + vm + wx)} \right|
			  \leq \frac{1}{2\pi} \left| \frac{2\pi x}{C} \cdot \frac{w}{\alpha + vm + wx} \right|
			  \leq \frac{1}{C} < 1.
\]
When $ \sigma > 0 $, the function $ g(\xi) $ is decreasing, and hence Lemma \ref{exp sum}
can be applied. 
For sufficiently large $ N $, we can take $ \varepsilon $ such that 
$ c -\varepsilon < 0 < d + \varepsilon < 1 $, by which only the term with $ \nu = 0 $ appears in the
sum on the right-hand side of (\ref{exponential-sum}). We obtain from (\ref{exponential-sum}) that
\[
	\sum_{x < n \leq N} \frac{e^{it\log(\alpha + vm + wn)}}{(\alpha + vm + wn)^{-\sigma}}
	= \int_x^N (\alpha + vm + w\xi)^{-\sigma + it} d\xi + O((m+x)^{-\sigma}).
\]
Taking complex conjugates on the both sides, we have
\begin{eqnarray*}
	\sum_{x < n \leq N} \frac{1}{(\alpha + vm + wn)^s}
	&=& \int_x^N (\alpha + vm + w\xi)^{-s} d\xi + O((m+x)^{-\sigma})	\\
	&=& \frac{(\alpha + vm + wN)^{1-s}-(\alpha + vm + wx)^{1-s}}{w(1-s)} + O((m+x)^{-\sigma}).
\end{eqnarray*}
Therefore, we obtain 
\begin{eqnarray*}
	B_1 &=& \sum_{m \leq x} \sum_{x < n \leq N} \frac{1}{(\alpha + vm + wn)^s}	\\
		&=& \sum_{m \leq x} \left\{ \frac{(\alpha + vm + wN)^{1-s}-(\alpha + vm + wx)^{1-s}}{w(1-s)} 
					+ O((m+x)^{-\sigma}) \right\}	\\
		&=& \frac{1}{w(1-s)} \left\{ \sum_{m \leq x} \frac{1}{(\alpha + vm + wN)^{s-1}}
								- \sum_{m \leq x} \frac{1}{(\alpha + vm + wx)^{s-1}}\right\}
					+ O(x^{1-\sigma}).				
\end{eqnarray*}
We denote the first and the second term on the right-hand side by 
$ (1/{w(1-s)})(B_{11}-B_{12}) $, and apply Lemma \ref{exp sum}
for $ B_{11} $ and $ B_{12} $. For $ B_{11} $ set
\[
	f(\xi) = \frac{t}{2\pi}\log(\alpha + v\xi + wN), \ g(\xi) = (\alpha + v\xi + wN)^{1-\sigma}
\]
and on taking $ (a,b) = (0,x) $ in Lemma \ref{exp sum}.
We can treat $ B_{12} $ similarly, where Lemma \ref{exp sum} is applied on replacing 
the variable $ \xi $ by $ \eta $, on setting
\[
	f(\eta) = \frac{t}{2\pi}\log(\alpha + v\eta + wx), \ g(x) = (\alpha + v\eta + wx)^{1-\sigma}
\]
and $ (a,b) = (0,x) $. Then we have
\begin{eqnarray*}
	&& B_{11} = \frac{(\alpha+vx+wN)^{2-s}-(\alpha+wN)^{2-s}}{v(2-s)} + O(N^{1-\sigma}),	\\
	&& B_{12} = \frac{(\alpha+vx+wx)^{2-s}-(\alpha+wx)^{2-s}}{v(2-s)} + O(x^{1-\sigma}).
\end{eqnarray*}
Therefore, we obtain
\begin{eqnarray*}
	&& B_1 = \frac{(\alpha+vx+wN)^{2-s}-(\alpha+vx+wx)^{2-s}-(\alpha+wN)^{2-s}+(\alpha+wx)^{2-s}}
{vw(s-1)(s-2)} \\
	&& \qquad \quad + O(N^{1-\sigma}) + O(x^{1-\sigma}).	
\end{eqnarray*}
By the argument similar to the treatment of $ B_1 $, we obtain
\begin{eqnarray*}
	&& B_2 = \frac{(\alpha+vN+wx)^{2-s}-(\alpha+vx+wx)^{2-s}-(\alpha+vN)^{2-s}+(\alpha+vx)^{2-s}}
{vw(s-1)(s-2)} \\
	&& \qquad \quad + O(N^{1-\sigma}) + O(x^{1-\sigma})
\end{eqnarray*}
and
\begin{eqnarray*}
	B_3 &=& \frac{(\alpha+vN+wN)^{2-s}-(\alpha+vx+wN)^{2-s}-(\alpha+vN+wx)^{2-s}+(\alpha+vx+wx)^{2-s}}{vw(s-1)(s-2)} \\
	&& \qquad + O(N^{1-\sigma}) + O(x^{1-\sigma}).	
\end{eqnarray*}
Summing up the results above, we obtain
\begin{eqnarray*}
	&& \sum_{m=0}^N \sum_{n=0}^N \frac{1}{(\alpha + vm + wn)^s}	\\
	&& \qquad = \sum_{m \leq x} \sum_{n \leq x} \frac{1}{(\alpha + vm + wn)^s} 
		+ \frac{(\alpha+vx)^{2-s}+(\alpha+wx)^{2-s} - (\alpha+vx+wx)^{2-s}}{vw(s-1)(s-2)}	\\
	&& \qquad - \frac{(\alpha+vN)^{2-s}+(\alpha+wN)^{2-s} - (\alpha+vN+wN)^{2-s}}{vw(s-1)(s-2)} 
		+ O(x^{1-\sigma}) + O(N^{1-\sigma}),
\end{eqnarray*}
and by (\ref{eq:sum_{x<N}}), we conclude that
\begin{eqnarray*}
	&& \zeta_2(s,\alpha;v,w)	\\
	&& \quad = \sum_{m \leq x} \sum_{n \leq x} \frac{1}{(\alpha + vm + wn)^s}
		+ \frac{(\alpha+vx)^{2-s}+(\alpha+wx)^{2-s} - (\alpha+vx+wx)^{2-s}}{vw(s-1)(s-2)}	\\
	&& \quad \ \ \ + O(x^{1-\sigma}) + O(N^{1-\sigma})
\end{eqnarray*}
in the region $ \sigma > 1 $. Letting $ N \rightarrow \infty $, we obtain the proof of 
Theorem \ref{order_of_Barnes_zeta}.
\qed

%%%%%%%% §4. Proof of Theorem 2 %%%%%%%%%%%%%%%%%%%%%%%%%%%%%%%%%%%%%%%%%%%
%%%%%%%%%%%%%%%%%%%%%%%%%%%%%%%%%%%%%%%%%%%%%%%%%%%%%%%%%%%%%%%%%%%%%%%%%%%%
%%%%%%%%%%%%%%%%%%%%%%%%%%%%%%%%%%%%%%%%%%%%%%%%%%%%%%%%%%%%%%%%%%%%%%%%%%%%
%%%%%%%%%%%%%%%%%%%%%%%%%%%%%%%%%%%%%%%%%%%%%%%%%%%%%%%%%%%%%%%%%%%%%%%%%%%%
%%%%%%%%%%%%%%%%%%%%%%%%%%%%%%%%%%%%%%%%%%%%%%%%%%%%%%%%%%%%%%%%%%%%%%%%%%%%
\section{Proof of Theorem 2}

In this section, we prove Theorem \ref{th:Main_Theorem2} from Theorem \ref{order_of_Barnes_zeta}.
\medskip

\textbf{Proof of Theorem \ref{th:Main_Theorem2}}. 
Setting $ C = 2\pi $ and $ x = t $ in (\ref{app_of_Barnes_zeta}), we easily see that
the second term on the right-hand side is $ O(t^{-\sigma}) $, hence we have
\begin{equation}\label{zeta_2 appro}
	\zeta_2(s, \alpha; v,w) = \sum_{m \leq t} \sum_{n \leq t} \frac{1}{(\alpha + vm + wn)^s}
							+ O(t^{1-\sigma}).
\end{equation}
We denote the first term on the right-hand side by $ \Sigma(s) $. Then
\begin{eqnarray*}
	&& \int_1^T |\Sigma(s)|^2 dt	\\
	&& \quad 
	= \int_1^T \sum_{m_1 \leq t} \sum_{n_1 \leq t} \frac{1}{(\alpha+ vm_1 + wn_1)^{\sigma +it}} 
			\sum_{m_2 \leq t} \sum_{n_2 \leq t} \frac{1}{(\alpha + vm_2 + wn_2)^{\sigma -it}} dt.
\end{eqnarray*}
Now we change the order of summation and integration. First we note that 
$ 1 \leq m_1, n_1, m_2, n_2 \leq T $.
Let us fix one such $ (m_1, n_1, m_2, n_2) $. Then from the condition 
$ m_1 \leq t,\ n_1 \leq t,\ m_2 \leq t,\ n_2 \leq t $, we find that the range of $ t $ is
$ M = \max\{m_1,n_1,m_2,n_2\} \leq t \leq T $. Therefore
\begin{eqnarray*}
	&& \int_1^T |\Sigma(s)|^2 dt	 \\
	&& \qquad = \sum_{m_1, n_1 \leq T} \frac{1}{(\alpha+ vm_1 + wn_1)^\sigma} 
			\sum_{m_2, n_2 \leq T} \frac{1}{(\alpha + vm_2 + wn_2)^\sigma}	\\
	&& \qquad \qquad \qquad \qquad \qquad \qquad \qquad \qquad \qquad \qquad \qquad \quad
		\times \int_M^T \left( \frac{\alpha + vm_1 + wn_1}{\alpha + vm_2 + wn_2} \right)^{it} dt \\
	&& \qquad = \mathop{\sum_{0 \leq m_1, n_1, m_2, n_2 \leq T}}
			\limits_{vm_1+wn_1 = vm_2+wn_2} 
			\frac{1}{(\alpha + v m_1 + w n_1)^\sigma(\alpha + v m_2 + w n_2)^\sigma} \times (T-M) \\
	&& \qquad + \mathop{\sum_{m_1, n_1, m_2, n_2 \geq 0}}
			\limits_{vm_1+wn_1 \neq vm_2+wn_2} 
			\frac{1}{(\alpha + v m_1 + w n_1)^\sigma(\alpha + v m_2 + w n_2)^\sigma} \\
	&& \qquad
			\times
			\frac{e^{iT\log\{(\alpha + v m_2 + w n_2)/(\alpha + v m_1 + w n_1)\}}
				- e^{iM\log\{(\alpha + v m_2 + w n_2)/(\alpha + v m_1 + w n_1)\}}}
			{i \log\{(\alpha + v m_2 + w n_2)/(\alpha + v m_1 + w n_1)\}}.	
\end{eqnarray*}
We denote the first and the second term on the right-hand side by $ S_1 $ and $ S_2 $ respectively.
As for $ S_1 $, we have
\begin{eqnarray*}
	S_1 &=& T \left\{ 
			\zeta_2^{[2]}(\sigma, \sigma, \alpha, v, w) 
			- ( U_1 + U_2 + U_3 + U_4 ) \right\}		\\
		&&  \qquad	
			- \mathop{\sum_{0 \leq m_1, n_1, m_2, n_2 \leq T}}
			\limits_{vm_1+wn_1 = vm_2+wn_2} 
			\frac{M(m_1,n_1,m_2,n_2)}{(\alpha + v m_1 + w n_1)^\sigma(\alpha + v m_2 + w n_2)^\sigma},
\end{eqnarray*}
where 
\begin{eqnarray*}
	U_1 &=& \left( 
			\sum_{\substack{m_1 > T \\ 0 \leq n_1, m_2, n_2 \leq T \\ vm_1+wn_1 = vm_2+wn_2}}
			+ \sum_{\substack{n_1 > T \\ 0 \leq m_1, m_2, n_2 \leq T \\ vm_1+wn_1 = vm_2+wn_2}}
		  \right.	\\
		& & \left.
			+ \sum_{\substack{m_2 > T \\ 0 \leq m_1, n_1, n_2 \leq T \\ vm_1+wn_1 = vm_2+wn_2}}
			+ \sum_{\substack{n_2 > T \\ 0 \leq m_1, n_1, m_2 \leq T \\ vm_1+wn_1 = vm_2+wn_2}}
			\right)
			\frac{1}{(\alpha + v m_1 + w n_1)^\sigma(\alpha + v m_2 + w n_2)^\sigma},	\\
	U_2 &=& \left( 
			\sum_{\substack{m_1, m_2 > T \\ 0 \leq n_1, n_2 \leq T \\ vm_1+wn_1 = vm_2+wn_2}}
			+ \sum_{\substack{n_1, n_2 > T \\ 0 \leq m_1, m_2 \leq T \\ vm_1+wn_1 = vm_2+wn_2}}
			\right)
			\frac{1}{(\alpha + v m_1 + w n_1)^\sigma(\alpha + v m_2 + w n_2)^\sigma},
\end{eqnarray*}
\begin{eqnarray*}
	U_3 &=& \left( 
			\sum_{\substack{m_1, n_2 > T \\ 0 \leq n_1, m_2 \leq T \\ vm_1+wn_1 = vm_2+wn_2}}
			+ \sum_{\substack{n_1, m_2 > T \\ 0 \leq m_1, n_2 \leq T \\ vm_1+wn_1 = vm_2+wn_2}}
			\right)
			\frac{1}{(\alpha + v m_1 + w n_1)^\sigma(\alpha + v m_2 + w n_2)^\sigma},	\\
	U_4 &=& \left( 
			\sum_{\substack{m_1, m_2, n_1 > T \\ 0 \leq n_2 \leq T \\ vm_1+wn_1 = vm_2+wn_2}}
			+ \sum_{\substack{m_1, m_2, n_2 > T \\ 0 \leq n_1 \leq T \\ vm_1+wn_1 = vm_2+wn_2}}
		  \right.	\\
		& & \left.
			+ \sum_{\substack{m_1, n_1, n_2 > T \\ 0 \leq m_2 \leq T \\ vm_1+wn_1 = vm_2+wn_2}}
			+ \sum_{\substack{m_2, n_1, n_2 > T \\ 0 \leq m_1 \leq T \\ vm_1+wn_1 = vm_2+wn_2}}
			\right)
			\frac{1}{(\alpha + v m_1 + w n_1)^\sigma(\alpha + v m_2 + w n_2)^\sigma}.
\end{eqnarray*}
We can estimate $ U_1 $ as follows. Since $ \alpha + v m + w n \asymp 1 + m + n $ we have
\begin{eqnarray*}
	U_1 \ll \sum_{\substack{{k>T} \\ 0 \leq l,m,n \leq T \\ k+l\asymp m+n}}\frac{1}{(1+k+l)^\sigma (1+m+n)^\sigma}.
\end{eqnarray*}
Setting $ j = k + l $, since $ T+1 < j \asymp m+n \leq 2T \ll T $ and $ m \ll j $, we obtain
\begin{eqnarray*}
	U_1 &\ll& \sum_{T+1 \leq j \ll T} \sum_{T < k \leq j } \sum_{0 \leq m \ll j} \frac{1}{(1+j)^{2\sigma}}	\\
		&\ll& \sum_{T+1 \leq j \ll T} \sum_{T<j \leq j} \frac{1}{(1+j)^{2\sigma -1}}
		=	  \sum_{{T+1 \leq j \ll T}} \frac{j-T}{(1+j)^{2 \sigma -1}}	\\
		&\ll& \sum_{{T+1 \leq j \ll T}} \frac{1}{(1+j)^{2 \sigma - 2}} 
		\ \asymp  \ T^{3-2 \sigma} \quad \left(\sigma \geq \frac{3}{2} \right).
\end{eqnarray*}
Similarly we obtain $ U_2,\ U_3,\ U_4 \ll T^{3-2\sigma} $. 
The sum involving $ M(m_1, n_1, m_2, n_2) $ in $ S_1 $ is estimated, since $ M(k, l, m, n) \ll k+l+m+n $ 
in this case, as
\begin{eqnarray*}
	& \ll & \sum_{\substack{0 \leq k, l, m, n \leq T \\ k + l \asymp m + n}}
		\frac{k+l+m+n}{(1+k+l)^{\sigma} (1+m+n)^{\sigma}}
		\ll  \sum_{0 \leq j \leq 2T} \sum_{0 \leq k \leq j} 
					\sum_{0 \leq m \ll j}\frac{1}{(1+j)^{2\sigma -1}}	\\
	& \ll &  \sum_{0 \leq j \leq 2T} \frac{1}{(1+j)^{2\sigma -3}} 
	\ll \begin{cases}
			T^{4-2\sigma}	 	& \left( 3/2 \leq \sigma < 2 \right),	\\
			\log{T}				& (\sigma = 2).
		\end{cases}	
\end{eqnarray*}
Therefore, we have
\begin{equation}\label{eq:S_1}
	S_1 = \zeta_2^{[2]}(\sigma,\sigma,\alpha,v,w)T + 
	\begin{cases}
			O(T^{4-2\sigma}) & 
				\left( 3/2 \leq \sigma < 2\right),	\\
			O(\log{T})	& (\sigma = 2).
	\end{cases}	
\end{equation}
Next, as for $ S_2 $, we have
\begin{eqnarray*}
	S_2 
	&\ll& \mathop{\sum_{0 \leq m_1, n_1, m_2, n_2 \leq T}}
			\limits_{vm_1+wn_1 < vm_2+wn_2}
			\frac{1}
				{(\alpha + v m_1 + w n_1)^\sigma(\alpha + v m_2 + w n_2)^\sigma}	\\
	&& \qquad \qquad \qquad \qquad \qquad \qquad 
		\times \frac{1}{\log\{(\alpha + v m_2 + w n_2)/(\alpha + v m_1 + w n_1)\}} 	\\
	&=& \left( \mathop{\sum_{0 \leq m_1, n_1, m_2, n_2 \leq T}}
			\limits_{\alpha + vm_1 + wn_1 < \alpha + vm_2 + wn_2 < 2(\alpha + vm_1 + wn_1)}
		+ \mathop{\sum_{0 \leq m_1, n_1, m_2, n_2 \leq T}}
			\limits_{ \alpha + vm_2 + wn_2  \geq 2(\alpha + vm_1 + wn_1)} 	
		\right) \\
	&&	\frac{1}{(\alpha + v m_1 + w n_1)^\sigma(\alpha + v m_2 + w n_2)^\sigma
				\log\{(\alpha + v m_2 + w n_2)/(\alpha + v m_1 + w n_1)\} }.
\end{eqnarray*}
We denote the first and the second term on the right-hand side by $ W_1 $ and $ W_2 $, respectively.
As for $ W_2 $, we have
\begin{eqnarray*}
	W_2 & \ll &  \mathop{\sum_{0 \leq m_1, n_1, m_2, n_2 \leq T}}
				\limits_{ \alpha + vm_2 + wn_2 \geq 2(\alpha + vm_1 + wn_1)}
				\frac{1}{(\alpha + v m_1 + w n_1)^\sigma(\alpha + v m_2 + w n_2)^\sigma}	 \\
		& \ll &  \left\{ \sum_{0 \leq m, n \leq T} \frac{1}{(1+m+n)^\sigma} \right\}^2		\\
		& = & \left\{ \sum_{0 \leq j \leq 2T} \sum_{m=0}^j \frac{1}{(1 + j)^\sigma} \right\}^2
		  =  \left\{ \sum_{0 \leq j \leq 2T} \frac{1}{(1 + j)^{\sigma -1}} \right\}^2		\\
		& \ll &	
			\begin{cases}
				T^{4-2\sigma}	& \left( 1 < \sigma < 2 \right),	\\
				(\log{T})^2		& (\sigma = 2).
			\end{cases}
\end{eqnarray*}
Next we consider the order of $ W_1 $.The range of $ n_2 $ in the inequalities 
$ \alpha + vm_1 + wn_1 < \alpha + vm_2 + wn_2 < 2(\alpha + vm_1 + wn_1) $
of the summation condition on $ W_1 $ is
\[
	\frac{v}{w}(m_1 - m_2) + n_1 < n_2 < 
	\min \left\{ \frac{\alpha}{w} + \frac{v}{w}(2m_1-m_2) + 2n_1,\ [T]+1 \right\}.
\]
Let $ \varepsilon = \varepsilon(m_1,m_2,n_1),  \delta = \delta (m_1,m_2,n_1) $ be
constants satisfying $ 0 \leq \varepsilon, \delta < 1 $ and
$ (v/w)(m_1-m_2) + n_1 + \varepsilon \in \Z $ and 
$ \alpha/w + (v/w)(2m_1-m_2) + 2n_1 - \delta \in \Z $.	 
Then
\[
	K = \min \left\{ \frac{\alpha}{w} + \frac{v}{w}\cdot m_1 + n_1 - \varepsilon - \delta ,\ 
			[T] - \frac{v}{w}(m_1 - m_2) - n_1 - \varepsilon  \right\}
\]
is an integer, and $ n_2 $ can be rewritten as
\[
	n_2 = \frac{v}{w}(m_1 - m_2) + n_1 + \varepsilon + k 
	\quad (\mathrm{for\ some\ } k=0,\ 1,\ 2,\ \ldots,\ K).
\] 
Since
\begin{eqnarray*}
	\log{\frac{\alpha + v m_2 + w n_2}{\alpha + v m_1 + w n_1}}
	= 	\log\left(1 + \frac{wk + w \varepsilon }{\alpha + v m_1 + w n_1}\right)
	\asymp  \frac{wk + w \varepsilon }{\alpha + v m_1 + w n_1},
\end{eqnarray*}
we obtain
\begin{eqnarray*}
	W_1 &\ll & \sum_{0 \leq m_1 \leq T} \sum_{0 \leq n_1 \leq T} 
			\sum_{0 \leq m_2 \ll K} \sum_{0 \leq k \leq K}
			\frac{1}{(\alpha + vm_1 +wn_1)^\sigma}	\\
		&&	\qquad \qquad \times \frac{1}{(\alpha + vm_1 +wn_1 + wk + w \varepsilon )^\sigma}
						\times \frac{\alpha + v m_1 + w n_1}{wk + w \varepsilon }		\\
		&\ll& \sum_{0 \leq m_1 \leq T} \sum_{0 \leq n_1 \leq T} \sum_{0 \leq m_2 \ll K} 
			\frac{\log{K}}{(\alpha + vm_1 + wn_1)^{2\sigma -1}}							\\
		&\ll& \sum_{0 \leq m_1 \leq T} \sum_{0 \leq n_1 \leq T} 
			\frac{\log{(\alpha + vm_1 + wn_1)}}{(\alpha + vm_1 + wn_1)^{2\sigma -2}}	\\
		&\ll& \sum_{0 \leq m \leq T} \sum_{0 \leq n \leq T} \frac{\log{(2+m+n)}}{(1+m+n)^{2\sigma -2}} \\
		& \ll &  \int_0^T \int_0^T \frac{\log{(2+x+y)}}{(1+x+y)^{2 \sigma - 2}} dxdy	\\
		& = &  
			\begin{cases}
				O(T^{4-2\sigma }\log{T})	& 
						\left( 1 < \sigma < 3/2, 3/2 < \sigma <  2 \right),	\\
				O(T(\log{T})^2)		& \left(\sigma = 3/2 \right),						\\
				O((\log{T})^2) 		& (\sigma = 2).
			\end{cases}
\end{eqnarray*}
Then, we have
\begin{eqnarray*}
	S_2 &=& \begin{cases}
				O(T^{4-2\sigma }) & \left( 1 < \sigma <  2 \right)	\\
				O((\log{T})^2) 		& (\sigma = 2)
			\end{cases}
			+
			\begin{cases}
				O(T^{4-2\sigma }\log{T})	& 
					\left( 1 < \sigma < 3/2, 3/2 < \sigma <  2 \right)	\\
				O(T(\log{T})^2)		& \left(\sigma = 3/2 \right)		\\
				O((\log{T})^2) 		& (\sigma = 2)
			\end{cases}													
			\\
		&=& \begin{cases}
				O(T^{4-2\sigma }\log{T})	& 
					\left( 1 < \sigma < 3/2, 3/2 < \sigma <  2 \right),	\\
				O(T(\log{T})^2)		& \left(\sigma = 3/2 \right),		\\
				O((\log{T})^2) 		& (\sigma = 2).
			\end{cases}
\end{eqnarray*}
By (\ref{eq:S_1}), we have
\begin{eqnarray*}
	&& \int_1^T |\Sigma(s)|^2 dt	\\
	&&  \qquad \quad 
	 = \zeta_2^{[2]}(\sigma,\sigma;\alpha,v,w)T +
	 \begin{cases}
				O(T^{4-2\sigma }\log{T})	& 
					\left( 1 < \sigma < 3/2, 3/2 < \sigma <  2 \right),	\\
				O(T(\log{T})^2)		& \left(\sigma = 3/2 \right),		\\
				O((\log{T})^2) 		& (\sigma = 2).
	 \end{cases}
\end{eqnarray*}
Furthermore, we obtain from (\ref{zeta_2 appro}) that
\begin{eqnarray*}
	&& \int_1^T |\zeta_2(\sigma + it,\alpha;v,w)|^2 dt 	\nonumber \\
	&& = \int_1^T|\Sigma(s) + O(t^{1-\sigma})|^2 dt	\nonumber \\
	&& = \int_1^T |\Sigma(s)|^2 dt + O \left( \int_1^T|\Sigma(s)|t^{1-\sigma} dt \right)
		+ O \left( \int_1^T t^{2-2\sigma} dt \right).	% \label{int|zeta|^2}
\end{eqnarray*}
We see that the third term on the right-hand side is estimated as
\begin{equation}
	\begin{cases}
			O(T^{3-2\sigma }) & ( 1 < \sigma < 3/2 ),	\\
			O(\log{T}) 		& (\sigma = 3/2),	\\
			O(1)			& (3/2 < \sigma \leq 2).
	\end{cases}
	\label{int|zeta|^2,2}
\end{equation}
Also, using the Cauchy-Schwarz inequality for the second term on the right-hand side, we see that
\begin{eqnarray}
	\int_1^T|\Sigma(s)|t^{1-\sigma} dt
	& \leq & \left( \int_1^T |\Sigma(s)|^2 dt \right)^{1/2} 
		\left( \int_1^T t^{2-2\sigma} dt \right)^{1/2}		\nonumber\\
	& = &
	\begin{cases}
				O(T^{7/2 - 2\sigma }(\log{T})^{1/2}) & \left( 1 < \sigma < 3/2 \right),	\\
				O(T^{1/2}(\log{T})^{3/2})		& \left(\sigma = 3/2 \right),		\\
				O(T^{1/2}) 		& (3/2 < \sigma \leq 2).
	 \end{cases}
	 \label{Cauchy-Schwarz}
\end{eqnarray}
Therefore, since (\ref{int|zeta|^2,2}) and (\ref{Cauchy-Schwarz}) we have
\begin{equation*}
	\int_1^T | \zeta_2(s, \alpha ; v,w)|^2 dt 
	= \zeta_2^{[2]}(\sigma, \sigma, \alpha ; v,w)T 
	+ 
	\begin{cases}
			O(T^{4-2\sigma} \log{T}) & 
				\left( 3/2 < \sigma \leq 7/4 \right),	\\
			O(T^{1/2})	& (7/4 < \sigma \leq 2),
	\end{cases}
\end{equation*}
and hence the proof of Theorem \ref{th:Main_Theorem2} is complete.
\qed

\bigskip \bigskip
\author{Takashi Miyagawa}:	\\
Graduate School of Mathematics, \\
Nagoya University, \\
Chikusa-ku, Nagoya, 464-8602 Japan	\\
E-mail: d15001n@math.nagoya-u.ac.jp

\bigskip

\begin{thebibliography}{99}
\bibitem{B1}
\textsc{E. W. Barnes}, The Genesis of the Double Gamma Functions, 
Proc. London Math. Soc. \textbf{31} (1899), 358--381.

\bibitem{B2}	
\textsc{E. W. Barnes}, The Theory of the Double Gamma Function, 
Philos. Trans. Roy. Soc. \textbf{197} (1901), 265--387.

\bibitem{B3}	
\textsc{E. W. Barnes}, On the Theory of the Multiple Gamma Function,
Trans. Cambridge Philos. Soc. \textbf{19} (1904), 374--425.

\bibitem{IMN}	
\textsc{S. Ikeda, K. Matsuoka, Y. Nagata}, 
On certain mean values of the double zeta-function, 
Nagoya Math. J. \textbf{217} (2015), 161--190.

\bibitem{KMT}
\textsc{Y. Komori, K. Matsumoto, H. Tsumura},
Barens multiple zeta-functions, Ramanujan's formula, and relevant series involving hyperbolic functions, 
J.Ramanujan Math. Soc. \textbf{28} (2013), 49--69.

\bibitem{MT}
\textsc{K. Matsumoto, H. Tsumura}, Mean value theorems for the double zeta-function,
J. Math. Soc. Japan \textbf{67} (2015), 383--408.

\bibitem{OO}
\textsc{T. Okamoto, T. Onozuka}, 
Mean value theorems for the Mordell-Tornheim double zeta-functions,
				Ramanujan J. \textbf{37} (2015), 131--163.

\bibitem{Tit}
\textsc{E. C. Titchmarsh}, 
The Theory of the Riemann Zeta-function. 2nd ed., Edited and
				with a preface by D. R. Heath-Brown, The Clarendon Press, Oxford University
				Press, New York, 1986.
\end{thebibliography}
\end{document}